\documentclass[a4paper,conference]{IEEEtran}
\IEEEoverridecommandlockouts

\usepackage{cite,color,url}
\usepackage{amsmath,amssymb,amsfonts,hyperref,multirow}
\usepackage{booktabs}
\usepackage{algorithmic,tikz}
\usepackage{subcaption,amsthm,longtable,lscape,amsmath}
\usepackage{graphicx,longtable}
\usepackage{textcomp}
\usepackage{xcolor}
\def\BibTeX{{\rm B\kern-.05em{\sc i\kern-.025em b}\kern-.08em
    T\kern-.1667em\lower.7ex\hbox{E}\kern-.125emX}}
\begin{document}

\title{On Solving the Minimum Spanning Tree \\Problem with
Conflicting Edge Pairs
}
\author{\IEEEauthorblockN{Roberto Montemanni}
\IEEEauthorblockA{\textit{Department of Sciences and Methods for Engineering} \\
\textit{University of Modena and Reggio Emilia}\\
Reggio Emilia, Italy \\
roberto.montemanni@unimore.it}
\and
\IEEEauthorblockN{Derek H. Smith}
\IEEEauthorblockA{\textit{Faculty of Computing, Engineering and Science} \\
\textit{University of South Wales}\\
Pontypridd, Wales, UK \\
derek.smith@southwales.ac.uk}}

\maketitle
\begin{abstract}
The Minimum Spanning Tree with  Conflicting Edge Pairs is a generalization that adds conflict constraints to a classical optimization problem on graphs used to model several real-world applications. In the last few years several approaches, both heuristic and exact, have been proposed to attack the problem. In this paper we consider a mixed integer linear program never approached before in the context of the problem under investigation, and we solve it with an open-source solver. Computational results on the benchmark instances commonly used in the literature of the problem are reported. The results indicate that the approach we propose, in its simplicity, obtains results aligned with those of the much more sophisticated approaches available. During the experimental campaign 6 instances have been closed for the first time, with 9 improved best-known lower bounds and 16 improved best-known upper bounds over the 230 instances considered.
\end{abstract}


\section{Introduction}
The Minimum Spanning Tree Problem with Conflicting Edge Pairs (MSTC) is
a very recent variant of the classical minimum spanning tree (MST) problem \cite{pri57}.
Given a connected undirected graph with costs on the edges and  a set of
edges pairs in conflict with each other, the MSTC consists of finding the spanning tree with minimum total cost that uses at most one of the edges of each conflicting pair. 
The MSTC was first introduced by Darmann et al. \cite{dar11}. It was shown that under general settings the problem cannot be solved in polynomial time, although some special cases are shown to be polynomially-solvable. Several heuristics for the problem were subsequently proposed in \cite{zha11}, while a branch-and-cut based on a mathematical programming model and some valid inequalities was presented in \cite{som14}. More recently, a metaheuristic method was designed and presented in \cite{car19}, while some new methods based on Lagrangian relaxations able to produce high-quality lower bounds for the optimal solution and heuristic solutions were introduced in \cite{car21}. Finally,  a new  branch-and-cut approach based on new valid inequalities was discussed in \cite{car19b}, and shown to represent the state-of-the-art for exact solving methods. 

MSTC can be used to represent real problem arising in different domains. For example,  the design of an offshore wind farm network, where the wind-turbines have to be connected together and the different possible point-to-point connections are characterized by different costs. Technical reasons prevent the use of overlapped cables in this domain \cite{kle15}, and this generates conflict constraints. Other practical applications are in the design and installation of long-distance pipelines, where conflicts typically represent the impossibility of using infrastructures of competitor service providers, or crossing countries with reciprocal political issues with a same line \cite{dar09}. 

In this paper, a compact mixed integer linear programming model for the problem is considered for the first time in the context of the MSTC and it is solved by the open-source solver CP-SAT, which is part of the Google OR-Tools \cite{cpsat} optimization suite. Successful application of this solver on optimization problems with  characteristics similar to the problem under investigation, motivated our study \cite{md23}, \cite{cor}, \cite{rm25}. A vast experimental campaign on the  benchmark instances previously used to validate algorithms in the literature is carried out.

The overall organization of the paper can be summarized as follows. Section \ref{desc} is dedicated to the formal introduction and definition of the Minimum Spanning Tree Problem with Conflicting Edge Pairs. Section \ref{model} presents a compact mixed integer linear programming model to represent the problem. In Section \ref{exp} the model proposed is analyzed from an experimental perspective. After the description of the benchmark instances considered and of the methods used for comparison purposes, the detailed results obtained by solving the new model are reported. Some conclusions about the work presented are drawn in Section \ref{conc} to close the paper.

\section{Problem Description}\label{desc}
The Minimum Spanning Tree Problem with Additional Conflicts
can be defined on an undirected graph $G = (V , E)$, with $V$ denoting
the set of nodes and $E$ the set of edges.  Each edge $\{i,j\} \in E$ is associated with a cost $u_{ij} \in \mathbb{N}^+$. Moreover, two edges
$\{i, j\}, \{k,l\} \in E$ are said to be \emph{in conflict} if at most one of them can be included in a feasible solution. The set of edges in conflict with the edge $\{i, j\} \in E$ is defined by the set $\delta(i, j)$. 
The objective of the MSTC is
to identify a spanning tree (a tree including all nodes)  fulfilling all  the conflict constraints with minimum possible total cost.

\begin{figure*}
{
\begin{center}
{
\begin{tikzpicture}[node distance={2.5cm}, main/.style = {draw, circle}]
			\node[main,minimum size=0.75cm] (f) {f};
			\node[main,minimum size=0.75cm] (a) [above right of=f] {a};
			\node[main,minimum size=0.75cm] (b) [below right of=f] {b};
			\node[main,minimum size=0.75cm] (c) [below right of=a] {c};
			\node[main,minimum size=0.75cm] (d) [above right of=c] {d};
			\node[main,minimum size=0.75cm] (e) [below right of=c] {e};
			\node[main,minimum size=0.75cm] (g) [below right of=d] {g};
			\draw [thick, line width=1.2,-] (f) to node[midway,above left] {3} (a) ;
			\draw [thick, line width=1.2,-,color=green] (f) to node[midway,below left] {6} (b) ;
			\draw [thick, line width=1.2,-,color=black] (a) to node[midway,below left] {1} (c) ;
			\draw [thick, line width=1.2,-,color=red] (a) to node[midway,above] {2} (d) ;
			\draw [thick, line width=1.2,-] (b) to node[midway,left] {1} (a) ;
			\draw [thick, line width=1.2,-,color=blue] (b) to node[midway,above left] {4} (c) ;
			\draw [thick, line width=1.2,-,color=blue] (b) to node[midway,below] {3} (e) ;
			\draw [thick, line width=1.2,-,color=red] (c) to node[midway,below right] {2} (d) ;
			\draw [thick, line width=1.2,-,color=black] (c) to node[midway,above right] {4} (e) ;
			\draw [thick, line width=1.2,-] (c) to node[midway,above] {3} (g) ;
			\draw [thick, line width=1.2,-] (d) to node[midway,above right] {5} (g) ;
			\draw [thick, line width=1.2,-,color=green] (e) to node[midway,below right] {7} (g) ;
\end{tikzpicture}
\hspace{1cm}
\begin{tikzpicture}[node distance={2.5cm}, main/.style = {draw, circle}]
			\node[main,minimum size=0.75cm] (s) {f};
			\node[main,minimum size=0.75cm] (a) [above right of=s] {a};
			\node[main,minimum size=0.75cm] (b) [below right of=s] {b};
			\node[main,minimum size=0.75cm] (c) [below right of=a] {c};
			\node[main,minimum size=0.75cm] (d) [above right of=c] {d};
			\node[main,minimum size=0.75cm] (e) [below right of=c] {e};
			\node[main,minimum size=0.75cm] (t) [below right of=d] {g};
			\draw [thick, line width=1.2,-] (s) to node[midway,above left] {3} (a) ;
			\draw [thick, line width=1.2,-,color=black] (a) to node[midway,below left] {1} (b) ;
			\draw [thick, line width=1.2,-,color=black] (a) to node[midway,below left] {1} (c) ;
			\draw [thick, line width=1.2,-] (c) to node[midway,above] {3} (t) ;
			\draw [thick, line width=1.2,-,color=red] (a) to node[midway,above right] {2} (d) ;
			\draw [thick, line width=1.2,-,color=white] (b) to node[midway,below] {3} (e) ;
			\draw [thick, line width=1.2,-,color=blue] (b) to node[midway,below] {3} (e) ;
\end{tikzpicture}
}
	\caption{On the left an example of a small MSTC instance is shown, where the edge capacities are placed by the edges. Edges in black are not affected by conflicts while any two edges sharing the same red, blue or green color are in conflict. On the right, an optimal solution that does not violate any conflict is shown. 			\vspace{-.5cm}
}
	\label{figu}
\end{center}
}
\end{figure*}
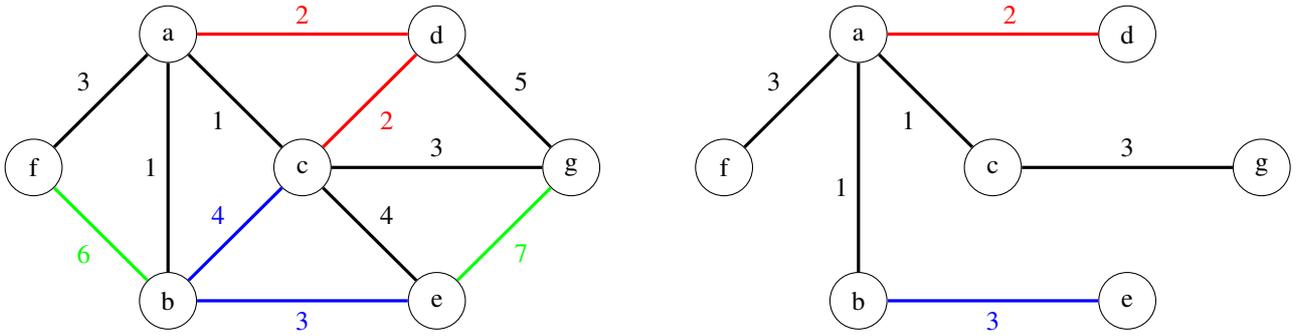

A small example of an MFPC instance and an optimal solution for the instance are provided in Figure \ref{figu}.

\section{A Mixed Integer Linear Programming Model}\label{model}
In this section a model for the MSTC, obtained by extending a classical model for the Minimum Spanning Tree Problem \cite{mag95}, is proposed. The model is based on multi-commodity flow and relies on a set $A$ of (directed) arcs, obtained by including the two arcs $(i,j)$ and ($j,i)$ for each edge $\{i,j\} \in E$. A node $s$ is arbitrarily selected as the root (node 0 in our implementation) and the aim of the model is to send one unit of flow from $s$ to any other node of $V$ using the arcs of $A$. Formally, a variable $x_{ij} \in \mathbb{Z}^+_0$ is introduced for each $(i,j) \in A$. It takes the value of flow transiting the arc $(i,j) \in A$. Another set of variables is required to match the multi-commodity flow problem back to the MSTC problem. Then a binary variable $y_{ij}$ is defined for each edge $\{i,j\} \in E$. It takes value 1 if the edge is part of the solution, 0 otherwise. The resulting compact model is as follows.
\begin{align} 
\ \ &  \min \ \  \sum_{\{i,j\} \in E} u_{ij}y_{ij}& \label{1}\\ 
s.t. \ \ &	\sum_{(j,i) \in A}  \!\!\!\! x_{ji}- \!\!\!\! \sum_{(i,j) \in A} \!\!\! x_{ij} \! = \! \begin{cases}
-|V|+1  \ \! \text{if $i=s$}\\
+1   \! \ \ \ \ \ \ \ \ \ \! \text{if $i \neq s$}
\end{cases} \!\!\!\!\!\!\!\!\!\!\!\!\!\!\!\!\!\!\!\!& i \in V\label{2}\\
& x_{ij} \leq (|V|-1) y_{ij} & \{i,j\} \in E \label{3a}\\
& x_{ji} \leq (|V|-1) y_{ij} & \{i,j\} \in E \label{3b}\\
& {y_{ij}+y_{kl} \leq 1} & \!\!\!\!\!\!\!\!\!\!\!\!\!\!\!\!\!\!\!\!\!\!\!\!\!\!\! \{i,j\} \in E, (k,l) \in \delta_{ij} \label{4}\\
& x_{ij} \ge 0 & (i,j) \in A \label{5}\\
& y_{ij} \in \{0,1\} & \{i,j\} \in E \label{6}
\end{align}
The objective function (\ref{1}) minimizes the cost of the spanning tree.
Constraints (\ref{2}) implement a classical multi-commodity flow model on the directed network, and impose that one unit of the commodity has to exit the source $s$ for each other node of $V$ (case $i=s$) and that each node different from $s$ has to retain exactly one unit of the commodity (case $i \neq s$). 
Inequalities (\ref{3a}) and (\ref{3b}) connect the $x$ and $y$ variables.
Constraints (\ref{4}) model conflicts and forbid conflicting edges to be active at the same time. 
Finally the domains of the variables are specified in constraints (\ref{5}) and (\ref{6}).

\section{Computational Experiments} \label{exp}
Section \ref{ben} is devoted to the description of the benchmark instances previously introduced in the literature, and used in the present study. An extensive comparison between the results achieved in this work and the results previously published is presented in Section \ref{res}.

\subsection{Benchmark Instances}\label{ben}
The computational tests are carried out on the benchmark sets originally proposed in \cite{zha11} and \cite{car19b}, that have been widely used in the previous literature to validate and compare solving methods. 

From the benchmark set proposed in \cite{zha11} we consider 50 challenging instances with a number of nodes ranging between 50 and 300 and a number of edges between 200 and 1000, with edge costs ranging between 0 and 500, and with different densities of conflicts.  The instances have been generated randomly.

The  newer set of benchmark instances introduced in  \cite{car19b} is formed by random instances with a number of nodes between 25 and 100, and  densities for the number of edges as 0.2, 0.3 or 0.4 (the previous set of benchmarks was characterized by lower densities). Costs are between 1 and 30 and the density of conflicts is either 0.01, 0.04 or 0.07. All the 180 instances of the set are considered for the experiments reported in the present work.

\subsection{Experimental Results} \label{res}
The model discussed in Section \ref{model} has been solved with the Google OR-Tools CP-SAT solver \cite{cpsat} version 9.12. The experiments have been run on a computer equipped with an Intel Core i7 12700F CPU, but the times have been normalized to the Intel Core i7-3770 CPU used for the experiments reported in \cite{car19}, \cite{car19b} and \cite{car21} and used as a reference for all the computation times reported. This allows a fairer comparison with the previous literature. The normalization of the computation times has been carried out according to the conversion ratio inferred from \url{http://gene.disi.unitn.it/test/cpu_list.php}. Note that a maximum computation time of 5010 seconds is allowed for all the methods for which a computation time is reported in the tables.

The first experiments have been carried out on the instances proposed in \cite{zha11}, and are summarized in Table \ref{t1}. The first columns identify each instance ($n$ is the number of nodes, $m$ is the number of edges and $p$ is the number of conflicts), while the rest of the table summarizes most of the results presented in the previous literature. The (not necessarily) optimal results obtained by solving a Mixed Integer Linear Programming (MILP) model and two lower bounding procedures (LB-MST and LB-MI) are taken from  \cite{zha11} (with no computation time). The lower and upper bounds (HDA and HDA+, respectively) from \cite{car21} are reported together with the upper bounds obtained with the metaheuristic method appeared in \cite{car19} and the lower and upper bounds obtained by the two branch-and-cut approaches discussed in \cite{som14} (without computation times) and \cite{car19b}, respectively. Finally, the results obtained by the method presented in this paper are presented in the last five columns (CP-SAT). On top of the upper and lower bounds obtained and the relative computation time, a percentage deviation (Dev \%) against the best lower and upper bounds is also reported, to better position the new approach. These deviations are calculated as $100 \cdot \frac{BK_{LB}-LB}{BK_{LB}}$ and $100 \cdot \frac{UB-BK_{UB}}{BK_{UB}}$ respectively, where $BK_{LB}$ and  $BK_{UB}$ are the best-known lower and upper bound available from the previous literature, and $LB$ and $UB$  are the values provided by CP-SAT. Negative deviations (improvements) are highlighted in bold.

The results of Table \ref{t1} suggests that the approach we propose (CP-SAT) is particularly effective in identifying infeasible instances: 5 of such instances were identified for the first time in this study by the powerful feasibility checker provided by the adopted solver (infeasible instances are marked as \emph{Infeas}) in the table). On top of this, another upper bound has been improved. In general, the new approach has performance comparable with those of the other state-of-the-art methods available, with gaps never above 8.5\% with respect to the best-known bounds. The average deviations reported are not very significant because heavily affected by the figures for the infeasibility proofs. Finally, it is worth observing how CP-SAT is able to provide meaningful lower and upper bounds for most of the instances, showing a  robustness that does not appear common with most of the other approaches.

A second set of experiments has been carried out on the instances proposed in \cite{car19b}, and the results are summarized in Table \ref{t2}, where the meaning of the column matches that of Table \ref{t1}, with the addition of parameter $s$ that is the seed used to generate the instances, and is used to distinguish among instances with the same characteristics. The results of all  the methods for which results are available on this dataset (the most recent ones) are reported.

The comparison proposed in Table \ref{t2} shows that the approach we propose (CP-SAT) performs well also on the instances proposed in \cite{car19b}, with an average deviation from the previously best-known results of 0.4\%  concerning the lower  bounds and 0.8\% concerning upper bounds. The maximum deviation is negligible for most of the instances, with only a few cases above 10\% concerning upper bounds and never above 2.2\% concerning lower bounds. A total of 4 new best-known lower bounds and 10 new best-known upper bounds have been produced during the experimental campaign, with one instance (\emph{50-245-2093-349}) closed for the first time.

In general, the new approach we propose appears to match the performance of state-of-the-art methods, notwithstanding its intrinsic simplicity, especially when compared to the  complexity of the ideas behind most of the other methods.

\section{Conclusions} \label{conc}
A compact mixed integer linear formulation for the Minimum Spanning Tree Problem with
Conflicting Edge Pairs was considered for the first time in the context of this problem, and solved via an open-source solver. 

The results of an experimental campaign, run on the benchmark instances adopted in the previous literature, are reported. The results indicate that the proposed approach, notwithstanding its intrinsic simplicity,  achieves computational results comparable with those of the substantially more complex methods proposed before,  with 6 instances closed for the first time, 9 best-known lower bounds and 16 best-known upper bounds improved, over the 230 instances considered.

\bibliographystyle{IEEEtran}
\bibliography{mybibfile}

\begin{thebibliography}{10}
\providecommand{\url}[1]{#1}
\csname url@samestyle\endcsname
\providecommand{\newblock}{\relax}
\providecommand{\bibinfo}[2]{#2}
\providecommand{\BIBentrySTDinterwordspacing}{\spaceskip=0pt\relax}
\providecommand{\BIBentryALTinterwordstretchfactor}{4}
\providecommand{\BIBentryALTinterwordspacing}{\spaceskip=\fontdimen2\font plus
\BIBentryALTinterwordstretchfactor\fontdimen3\font minus
  \fontdimen4\font\relax}
\providecommand{\BIBforeignlanguage}[2]{{%
\expandafter\ifx\csname l@#1\endcsname\relax
\typeout{** WARNING: IEEEtran.bst: No hyphenation pattern has been}%
\typeout{** loaded for the language `#1'. Using the pattern for}%
\typeout{** the default language instead.}%
\else
\language=\csname l@#1\endcsname
\fi
#2}}
\providecommand{\BIBdecl}{\relax}
\BIBdecl

\bibitem{pri57}
R.~C. Prim, ``Shortest connection networks and some generalizations,''
  \emph{The Bell System Technical Journal}, vol.~6, no.~36, pp. 1389--1401,
  1957.

\bibitem{dar11}
A.~Darmann, U.~Pferschy, and G.~Schauer, J.and~Woeginger, ``Paths, trees and
  matchings under disjunctive constraints,'' \emph{Discrete Applied
  Mathematics}, vol.~16, no. 159, pp. 1726--1735, 2011.

\bibitem{zha11}
R.~Zhang, S.~Kabadi, and A.~Punnen, ``The minimum spanning tree problem with
  conflict constraints and its variations,'' \emph{Discrete Optimization},
  vol.~2, no.~8, pp. 191--205, 2011.

\bibitem{som14}
P.~Samer and S.~Urrutia, ``A branch and cut algorithm for minimum spanning
  trees under conflict constraints,'' \emph{Optimization Letters}, vol.~1,
  no.~9, pp. 41--55, 2014.

\bibitem{car19}
F.~Carrabs, C.~Cerrone, and R.~Pentangelo, ``A multi-ethnic genetic approach
  for the minimum conflict weighted spanning tree problem,'' \emph{Networks},
  vol.~2, no.~74, pp. 134--147, 2019.

\bibitem{car21}
F.~Carrabs and M.~Gaudioso, ``A lagrangian approach for the minimum spanning
  tree problem with conflicting edge pairs,'' \emph{Networks}, vol.~1, no.~78,
  pp. 32--45, 2021.

\bibitem{car19b}
F.~Carrabs, R.~Cerulli, R.~Pentangelo, and A.~Raiconi, ``Minimum spanning tree
  with conflicting edge pairs: a branch-and-cut approach,'' \emph{Annals of
  Operations Research}, no. 298, pp. 65--78, 2019.

\bibitem{kle15}
A.~Klein, D.~Haugland, J.~Bauer, and M.~Mommer, ``An integer programming model
  for branching cable layouts in oﬀshore wind farms,'' \emph{Advances in
  Intelligent Systems and Computing}, no. 359, pp. 27--36, 2015.

\bibitem{dar09}
A.~Darmann, U.~Pferschy, and J.~Schauer, ``Determining a minimum spanning tree
  with disjunctive constraints,'' \emph{Lecture Notes in Computer Science}, no.
  5783, pp. 414--423, 2009.

\bibitem{cpsat}
L.~Perron and F.~Didier, ``{Google OR-Tools - CP-SAT},'' Google, 2025,
  \url{https://developers.google.com/optimization/cp/cp_solver/} [Accessed:
  2025-03-18].

\bibitem{md23}
R.~Montemanni and M.~Dell'Amico, ``Solving the parallel drone scheduling
  traveling salesman problem via constraint programming,'' \emph{Algorithms},
  vol.~16, no.~1, p.~40, 2023.

\bibitem{cor}
R.~Montemanni, M.~Dell'Amico, and A.~Corsini, ``Parallel drone scheduling
  vehicle routing problems with collective drones,'' \emph{Computers \&
  Operations Research}, vol. 163, p. 106514, 2024.

\bibitem{rm25}
R.~Montemanni, ``Solving a home healthcare routing and scheduling problem with
  real-world features,'' in \emph{Proceedings of the 9th International
  Conference on Machine Learning and Soft Computing}.\hskip 1em plus 0.5em
  minus 0.4em\relax Springer, to appear, 2025.

\bibitem{mag95}
T.~L. Magnanti and L.~Wolsey, \emph{Optimal trees}.\hskip 1em plus 0.5em minus
  0.4em\relax In \emph{Handbook in Operations Research and Management Science},
  M. O. Ball et al. (eds.), North Holland, Amsterdam, 1995, pp. 503--615.

\end{thebibliography}

\onecolumn
\begin{landscape}

\begin{table}[]
     \caption{Computational results for the instances introduced in\cite{zha11}.}\label{t1}
{ 
\scriptsize
\begin{tabular}{lrrrrrrrrrrrrrrrrrrrrr}
\toprule
\multicolumn{3}{c}{Instances} & \multicolumn{1}{c}{MILP} & \multicolumn{1}{c}{LB-MST} & \multicolumn{1}{c}{LB-MI} & \multicolumn{2}{c}{HAD} & \multicolumn{2}{c}{Mega} & \multicolumn{2}{c}{HAD+} & \multicolumn{2}{c}{BC} & \multicolumn{3}{c}{BC} & \multicolumn{5}{c}{CP-SAT}     \\
\multicolumn{3}{c}{} & \multicolumn{1}{c}{[5]} & \multicolumn{1}{c}{[5]} & \multicolumn{1}{c}{[5]} & \multicolumn{2}{c}{[2]} & \multicolumn{2}{c}{[4]} & \multicolumn{2}{c}{[2]} & \multicolumn{2}{c}{[3]} & \multicolumn{3}{c}{[1]} & \multicolumn{3}{c}{} &\multicolumn{2}{c}{Dev \%}     \\
n       & m       & p         & UB                              & LB                                 & LB                                & LB               & Sec          & UB               & Sec           & UB              & Sec            & LB              & UB           & LB        & UB      & Sec      & LB      & UB     & Sec     &   LB & UB \\
\cmidrule(lr){1-3}\cmidrule(lr){4-4}\cmidrule(lr){5-5}\cmidrule(lr){6-6}\cmidrule(lr){7-8}\cmidrule(lr){9-10}\cmidrule(lr){11-12}\cmidrule(lr){13-14}\cmidrule(lr){15-17}\cmidrule(lr){18-22}

50      & 200     & 199       & 708                              & 701.1                              & 702.8                             & 705.5            & 0.3          & 708              & 0.7           & 726             & 0.4            & 708.0           & 708          & 708.0     & 708     & 0.2      & 708     & 708    & 5.1     & 0.0                  & 0.0             \\
50      & 200     & 398       & 770                              & 739.8                              & 757.8                             & 761.0            & 0.3          & 770              & 0.7           & 780             & 0.5            & 770.0           & 770          & 770.0     & 770     & 0.5      & 770     & 770    & 9.2     & 0.0                  & 0.0             \\
\bottomrule
\end{tabular}
}
\end{table}

\end{landscape}

{\scriptsize
\begin{longtable}{lrrrrrrrrrrrrrrr}
     \caption{Computational results  for the instances introduced in \cite{car19b}.}\\  
     \toprule 
\multicolumn{4}{c}{Instances}                                                                            & \multicolumn{2}{c}{HDA } & \multicolumn{2}{c}{HDA+ } & \multicolumn{3}{c}{BC } & \multicolumn{5}{c}{CP-SAT}             \\ 
\multicolumn{4}{c}{}                                                                              & \multicolumn{2}{c}{{[}2{]}} & \multicolumn{2}{c}{ {[}2{]}} & \multicolumn{3}{c}{ {[}1{]}} & \multicolumn{3}{c}{}& \multicolumn{2}{c}{Dev \%}              \\ 
$n$&$m$& $p$& $s$    & LB              & Sec           & UB             & Sec             & LB          & UB       & Sec        & LB           & UB            & Sec &LB&UB    \\
\cmidrule(lr){1-4}\cmidrule(lr){5-6}\cmidrule(lr){7-8}\cmidrule(lr){9-11}\cmidrule(lr){12-16}
\endfirsthead
     \caption{Computational results   for the instances introduced in \cite{car19b} (continued).}\\  
     \toprule 
\multicolumn{4}{c}{Instances}                                                                            & \multicolumn{2}{c}{HDA } & \multicolumn{2}{c}{HDA+ } & \multicolumn{3}{c}{BC } & \multicolumn{5}{c}{CP-SAT}             \\ 
\multicolumn{4}{c}{}                                                                              & \multicolumn{2}{c}{{[}2{]}} & \multicolumn{2}{c}{ {[}2{]}} & \multicolumn{3}{c}{ {[}1{]}} & \multicolumn{3}{c}{}& \multicolumn{2}{c}{Dev \%}              \\ 
$n$&$m$& $p$& $s$    & LB              & Sec           & UB             & Sec             & LB          & UB       & Sec        & LB           & UB            & Sec &LB&UB    \\
\cmidrule(lr){1-4}\cmidrule(lr){5-6}\cmidrule(lr){7-8}\cmidrule(lr){9-11}\cmidrule(lr){12-16}
\endhead
\multicolumn{16}{r}{\textit{(Continued on next page)}} 
\endfoot
\endlastfoot
\label{t2}
 25           & 60            & 18     & 1    & 347.0  & 0.04   & 347 & 0.05 & 347.0  & 347  & 0.0    & 347          & 347  & 0.23    & 0.0           & 0.0             \\
25           & 60            & 18     & 7    & 389.0  & 0.04   & 389 & 0.04 & 389.0  & 389  & 0.0    & 389          & 389  & 0.18    & 0.0           & 0.0             \\
25           & 60            & 18     & 13   & 353.0  & 0.04   & 353 & 0.04 & 353.0  & 353  & 0.0    & 353          & 353  & 0.19    & 0.0           & 0.0             \\
25           & 60            & 18     & 19   & 346.0  & 0.04   & 346 & 0.05 & 346.0  & 346  & 0.0    & 346          & 346  & 0.19    & 0.0           & 0.0             \\
25           & 60            & 18     & 25   & 336.0  & 0.04   & 336 & 0.04 & 336.0  & 336  & 0.0    & 336          & 336  & 0.18    & 0.0           & 0.0             \\
25           & 90            & 71     & 31   & 379.6  & 0.09   & 381 & 0.07 & 381.0  & 381  & 0.0    & 381          & 381  & 0.45    & 0.0           & 0.0             \\
25           & 60            & 71     & 37   & 381.5  & 0.05   & 390 & 0.07 & 390.0  & 390  & 0.1    & 390          & 390  & 0.32    & 0.0           & 0.0             \\
25           & 60            & 71     & 43   & 372.0  & 0.05   & 372 & 0.06 & 372.0  & 372  & 0.0    & 372          & 372  & 0.44    & 0.0           & 0.0             \\
\bottomrule
\end{longtable}
}
\twocolumn

\end{document}